\documentstyle[12pt]{article} \vsize=28.7truecm \hsize=23truecm
 \def\vbar{\mathchoice{\vrule height2.3ptdepth-.3ptwidth.12pt\kern-
 .10pt}
    {\vrule height6.3ptdepth-.3ptwidth.11pt\kern-.11pt}
    {\vrule height5.1ptdepth-.30ptwidth.8pt\kern-.8pt}
    {\vrule height4.1ptdepth-.24ptwidth.6pt\kern-.7pt}}
\def\reel{\hbox{I\hskip-2pt R}}
\def\comp{\hbox{I\hskip-6pt C}}

\def\natu{\hbox{I\hskip-2pt N}}

\def\<{\langle}
\def\>{\rangle}

\def\n{{\boldmath n}}

\def\<{\langle}
\def\>{\rangle}
\def\reel{\hbox{I\hskip -2pt R}}

\def\n{{\noindent}}
\newtheorem{theorem}{Theorem}[section]

\newtheorem{prop}[theorem]{Proposition}
\newenvironment{Pff}{\hspace*{-\parindent}{\bf Proof}}{\hfill $\Box$
\vspace*{0.2cm}}

\begin{document}
 \vspace{2cm}

\title{{\bf
 {Mixture of the Riesz distribution with respect to the
 multivariate
 Poisson}}}
\author{Abdelhamid Hassairi\footnote{Corresponding author.
 \textit{E-mail address: abdelhamid.hassairi@fss.rnu.tn}}$\;$ and
Mahdi Louati
\\{\footnotesize{\it
Laboratory of Probability and Statistics. Sfax Faculty of Sciences,
B.P. 802, Tunisia.}}}
\date{}
\maketitle{Running title: \emph{Mixture of the Riesz distribution}}

\n $\overline{\hspace{15cm}}$\vskip0.3cm \n {\small {\bf Abstract}}
{\small The aim of this paper is to study the mixture of the Riesz
distribution on symmetric matrices with respect to the multivariate
Poisson distribution. We show, in particular, that this distribution
is related to the modified Bessel function of the first kind. We
also study the generated natural exponential family. We determine
the domain of the means and the variance function of this family.}
\vskip0.2cm\n{\small{\it{Keywords:}} Mixed distribution, Riesz
distribution, Bessel function, natural exponential family, variance function.\\
$\overline{\hspace{15cm}}$\vskip1cm
\section {Introduction}
Let $\mu_\lambda$ be a probability distribution on a finite
dimensional linear space $E$ depending on a parameter $\lambda$
which belongs to a subset $\Lambda$ of $\reel^r$. Suppose that
$$\mu_\lambda=f(x,\lambda)\sigma(dx),$$ where $\sigma$ is some
reference measure, and that for each $x$ in $E$, the map
$\lambda\mapsto f(x,\lambda)$ defined on $\Lambda$ is measurable.
For a probability distribution $\nu(d\lambda)$ on the set $\Lambda$,
define $$h(x)=\displaystyle\int_{\Lambda}
f(x,\lambda)\nu(d\lambda).$$ Then the probability measure
$$\mu_\nu(dx)=h(x)\sigma(dx)$$ is called the mixture of
the distribution $\mu_\lambda$ with respect to $\nu$. (See Feller
[4], Vol. II, page 53 or Johnson \textit{et al.} [8], page 360).
Usually, $\mu_\lambda$ is said the mixed distribution and $\nu$ the
mixing distribution (see Karlis and Meligkotsidou [9]).\vskip0.1cm\n
A special case of interest is when $\mu$ is not concentrated on an
affine hyperplane of $E$, $\mu_\lambda$ is the $\lambda-$power of
convolution of $\mu$, and $\Lambda$ is the so called Y{\o}rgensen
set of $\mu$. Specifically, let
\begin{equation}\label{MI1}
L_\mu(\theta)=\displaystyle\int_{E} \exp(\langle\theta,
x\rangle)\mu(dx)
\end{equation}
denote the Laplace transform of $\mu$, where $\langle,\rangle$ is
the duality crocket, and suppose that the set
\begin{equation}\label{MI2}
\Theta(\mu)=\textrm{interior}\{\theta \in E^\ast; \ L_\mu
(\theta)<+\infty\}
\end{equation}is nonempty. Then the set
\begin{equation}\label{MIN2}
\Lambda=\{\lambda>0; \ \exists \ \  \mu_\lambda \textrm{ such that }
L_{\mu_{_\lambda}}(\theta)=\left(L_{\mu}(\theta)\right)^\lambda
\textrm{ for all }\theta \in \Theta(\mu) \}
\end{equation} is called the Y{\o}rgensen set of $\mu$ and the measure $\mu_\lambda$ is its $\lambda-$power of
convolution. This parameter $\lambda$ appears in the most common
models, it is in particular, the intensity in a Poisson model, the
variance in a Gaussian model, and the shape parameter in a gamma
model. For $\lambda$ and $\lambda'$ in $\Lambda$, we have that
$$\mu_\lambda \ast\mu_{\lambda'}=\mu_{\lambda+\lambda'}.$$ The set
$\Lambda$ contains always the set $\natu^\ast$ of positive integers
and it is equal to $]0,+\infty[$ if and only if $\mu$ is infinitely
divisible (see Seshadri [13], page 155).\vskip0.1cm\noindent In
fact, if $X_1,\ldots,X_N$ are iid random variables with distribution
$\mu$, then the distribution of $X_1 +\ldots+X_N$ is $\mu_N
=\mu^{\star N}$, the $N-$power of convolution of $\mu$. Accordingly,
for any distribution $\mu$ and any positive integer $N$, one may
consider the distribution $\mu_N$ as defined in (\ref{MIN2}). When
$\mu$ is discrete, i.e., with countable support, then the mixture of
$\mu$ with respect of a distribution $\nu$ on the parameter $N$ is
known as a compound distribution. The most famous compound
distribution is the one corresponding to the case where $\nu$ is
Poisson (see Feller [4], Vol. I, page 286 or Vol. II, page 451 or
Aalen [1], or Y{\o}rgensen [14], page 140). In fact, the real
Poisson distribution appears in numerous works either as a mixed
distribution (see Bhattacharya and Holla [2] or Johnson \textit{et
al.} [8], page 366) or as a mixing distribution (see Perline [12]).
In the present work, we will be interested in a very special case in
which the mixed distribution is defined on the cone of $(r,r)$
positive definite symmetric matrices $\Omega$ with a parameter which
belongs to a subset of $\reel^r$. More precisely, the mixed
distribution will be the absolutely continuous Riesz model
introduced in Hassairi and Lajmi [6]
$$\left\{R(s,\sigma ), \ s\in
\displaystyle\prod_{i=1}^r\left]\frac{i-1}{2},+\infty\right[\right\}.$$
This model contains the Wishart model, since $R(s,\sigma )$ reduces
to a Wishart distribution when $s_1=s_2=\ldots=s_r$, and it has a
convolution property which is analogous to the one satisfied by the
ordinary powers of convolution. In fact, if $s$ and $s'$ are in
$\displaystyle\prod_{i=1}^r\left]\frac{i-1}{2},+\infty\right[$, then
$$R(s,\sigma)\ast R(s',\sigma)=R(s+s',\sigma).$$ The mixing distribution will be the
multivariate Poisson distribution on $\natu^r$. For simplicity, we
will be interested in the case where $\sigma$ is equal to the
identity matrix of size $r$ denoted $I_{r}$. We first show that the
mixture distribution is expressed in terms of the modified Bessel
function. We then determine the domain of the means of the generated
natural exponential family, and we calculate its variance function.
This provides a rich class of natural exponential families.
\section{The Riesz exponential dispersion model}
A dispersion model is a class of natural exponential families where
each family is a power of convolution of the other. In this section,
we will first recall some general facts concerning the exponential
dispersion models in an Euclidean space. Then we introduce the Riesz
dispersion model on symmetric matrices.
\subsection{Exponential dispersion model} Let $E$ be an
Euclidean space with finite dimension $n,$ and let $
\langle,\rangle$ denote the scalar product in $E$. If $\mu$ is a
positive measure on $E$, we denote by ${\mathcal{M}}(E)$ the set of
measures $\mu$ such that $\Theta(\mu) $ given in (\ref{MI2}) is not
empty and $\mu$ is not concentrated on an affine hyperplane of $E$.
If $\mu$ is in ${\mathcal{M}}(E)$, we denote \vskip0.2cm$\hfill
k_\mu (\theta)=\log L_\mu (\theta), \textrm{ for all }\theta
\textrm{ in }\Theta(\mu)\hfill$\vskip0.2cm\n the cumulant function
of $\mu$, where $L_\mu$ is the Laplace transform of $\mu$ defined in
(\ref{MI1}). \vskip0.1cm\n To each $\mu$ in ${\mathcal{M}}(E)$ and
$\theta$ in $\Theta(\mu)$, we associate the probability distribution
on $E$ \vskip0.2cm$\hfill
P(\theta,\mu)(dx)=exp\left(\langle\theta,x\rangle -k_\mu
(\theta)\right)\mu(dx).\hfill$\newline
 The set\vskip0.1cm$\hfill
 F=F(\mu)=\{P(\theta,\mu); \ \theta\in\Theta(\mu)\}\hfill$
\vskip0.3cm\n is called the natural exponential family (NEF)
generated by $\mu$. We also say that $\mu$ is a basis of $F$.
\vskip0.2cm\n The function $k_\mu$ is strictly convex and real
analytic. Its first derivative $k'_\mu$ defines a diffeomorphism
between $\Theta(\mu)$ and its image $M_F$. Since $k'_\mu
(\theta)=\displaystyle\int_E x P(\theta,\mu) (dx)$, $M_F$ is called
the domain of the means of $F$. The inverse function of $k'_\mu$ is
denoted by $\psi_\mu$ and setting $P(m,F)=P(\psi_\mu(m),\mu)$ the
probability of $F$ with mean $m$, we have $$F=\left\{P(m,F);m\in
M_F\right\},$$ which is the parametrization of $F$ by the
mean.\vskip0.1cm\n Now the covariance operator of $P(m,F)$ is
denoted by $V_F (m)$ and the map $$M_F\longrightarrow L_s(E); \
m\longmapsto V_F (m)=k''_\mu(\psi_\mu(m))$$ is called the variance
function of the NEF $F$. It is easy proved that for all $m\in M_F$,
$$V_F(m)=(\psi'_\mu(m))^{-1},$$ and an important feature of $V_F$ is
that it characterizes $F$ in the following sense: If $F$ and $F'$
are two NEFs such that $V_F (m)$ and $V_{F'}(m)$ coincide on a
nonempty open set of $M_F \cap M_{F'},$ then $F=F'$. In particular,
knowledge of the variance function gives knowledge of the natural
exponential family.\vskip0.1cm\noindent Let $\mu$ be an element of
${\mathcal{M}}(E)$ and let $ \Lambda$ be its Y{\o}rgensen set
defined by (\ref{MIN2}). Then the set
$$\{P(\theta,\lambda)= \exp\left(<\theta ,x>
-k_{\mu_\lambda}(\theta)\right) \mu; \ \theta\in\theta(\mu), \
\lambda\in \Lambda\}$$ is called the dispersion model generated by
$\mu$. For more details, we refer to Letac [11].
\subsection{Riesz natural exponential families}
Let $E$ be the Euclidean space of $(r,r)$ real symmetric matrices
equipped with the scalar product $\langle
x,y\rangle=\textrm{tr}(xy)$, and the inner product
$x.y=\displaystyle\frac{1}{2}\Big(xy+yx\Big),$ where $xy$ is the
ordinary product of two matrices.\vskip0.1cm\n We denote by
$e_1,e_2,\ldots,e_r$ the canonical basis of ${\reel}^r$;
$e_i=(0,\ldots,0,1,0\ldots0)$, (1 in the i$^\textrm{th}$ place), and
we set $c_i=diag(e_i)$ for all $1\leq i\leq r.$\vskip0.1cm\n For
$x\in E$, we consider the endomorphism $L(x)$ of $E$ defined by
$$L(x): \ y\longmapsto x.y$$ and we set
$$P(x)=2(L(x))^{^{2}}-L(x^{^{2}}).$$ We denote by
$\Omega$ the cone of $(r,r)$ real symmetric positive definite
matrices. \vskip0.1cm\n For $x= (x_{ij} )_{1\leq i,j \leq r}$ in $E$
and $1\leq k\leq r$, we define the sub-matrices
$$P_k(x)=(x_{ij} )_{1\leq i,j \leq k} \textrm{ and }
P_k^\ast(x)=(x_{ij} )_{r-k+1\leq i,j \leq r}.$$ For convenience,
$P_k(x)$ and $P_k^\ast(x)$ are also considered as elements of the
space $E$, where the other entries are equal to zero and we set
$P_0^\ast(x)=0.$\vskip0.1cm\n Let $\Delta_k(x)$ and
$\Delta^\ast_k(x)$ denote the determinant of $P_k(x)$ and the
determinant of $P^\ast_k(x)$ respectively. Then the generalized
power of $x$ in $\Omega$ is defined, for
$s=(s_1,s_2,\ldots,s_r)\in{\reel}^r$, by
\begin{equation}\label{FATM1}
\Delta_s(x)=\Delta_1(x)^{s_1-s_2}\Delta_2(x)^{s_2-s_3}\ldots\ldots\Delta_r(x)^{s_r}.
\end{equation}
Note that if for all $i\in\{1,\ldots,r\}$, $s_i=p\in {\reel}$, then
$\Delta_s(x)=(\det x)^p$. We also define
\begin{equation}\label{N1}
\Delta^\ast_s(x)=(\Delta^\ast_1(x))^{s_1-s_2}(\Delta^\ast_2(x))^{s_2-s_3}
\ldots(\Delta^\ast_{r-1}(x))^{s_{r-1}-s_r} (\Delta^\ast_r(x))^{s_r}.
\end{equation}
It is shown (see Hassairi and Lajmi [6]) that for all $x\in\Omega$
and all $s\in{\reel}^r$, we have
\begin{equation}\label{N0}
\Delta_s(x^{-1})=\Delta^\ast_ {-s^\ast}(x),
\end{equation}where $s^\ast=(s_r,s_{r-1},\ldots,s_1).$ .\vskip0.15cm\n We denote by ${\mathcal{T}}_{l}^+$ the set of lower
triangular matrices with positive diagonal elements, and for $u\in
{\mathcal{T}}_{l}^+,$ we define on $E$ the automorphism
$$u(y)=u yu^\ast,$$ where $u^\ast $ denotes the transpose matrix of
$u.$ \vskip0.1cm\n It is well known that for all $x\in\Omega,$ there
exists a unique $u\in {\mathcal{T}}_{l}^+$ such that
$$x=u(I_r),$$ where $I_r$ is the identity matrix of order $r,$ it
is the Cholesky decomposition of $x$.\vskip0.1cm\n We also have (see
Hassairi and Lajmi [6]) that for all $1\leq i\leq r$,
\begin{equation}\label{N15}
\left(P_i^{\ast}\left((u(I_r))^{-1}\right)\right)^{-1}
=u\left(\displaystyle\sum_{k=r-i+1}^r c_{k}\right),
\end{equation}
and for all $s=(s_1,s_2,\ldots,s_r)\in\reel^r,$
\begin{equation}\label{R1}
\Delta_ s(u(I_r))=\Delta^\ast_ {-s^\ast}(u^{\ast -1}(I_r)).
\end{equation}
\vskip0.1cm\noindent Recall also that for $x\in \Omega$ and $u\in
{\mathcal{T}}_{l}^+$, we have
\begin{equation}\label{Nour1}
\Delta_i(u(x))=\Delta_i(u(I_r))\Delta_i(x)=u_1^2\ldots
u_i^2\Delta_i(x),
\end{equation}where for all
$i\in\{1,\ldots,r\}$,
\begin{equation}\label{III1}
u_i=\langle u,c_i\rangle.
\end{equation}(See
Faraut and Koranyi [3], page 114).\vskip0.1cm\noindent
 Now let $\Xi$ be the set of elements $s=(s_1,s_2,\ldots,s_r)\in{\reel}^r$
defined as follows:\vskip0.3cm\noindent Consider the function $\xi$
defined from ${\reel}^+$ into $\{0;1\}$, by $$\xi: \ a\longmapsto
\left\{\begin{array}{cc}
% after \\: \hline or \cline{col1-col2} \cline{col3-col4} ...
  0 & \textrm{if }a=0, \\
 1 & \textrm{if }a>0. \\
\end{array}\right. $$
For all $(u_1,u_2,\ldots,u_r)\in{\reel}_+^{r},$ we define
$$\left\{
\begin{array}{ccc}
s_1=u_1 \\
           \\
s_k  = & u_k+\displaystyle\frac{\xi(u_1)+...+\xi(u_{k-1})}{2},
  \ \forall \ k\in\{2,\ldots,r\}. \\
\end{array}\right.$$
A result due to Gindikin [5] and proved in Faraut and Kor\'{a}nyi
[3], page 124, says that there exists a positive measure $R_s$ such
that for all $\theta\in -\Omega,$
$$L_{R_s}(\theta)=\int_E
e^{\langle\theta,x\rangle}R_s(dx)=\Delta_s(-\theta^{-1})$$ if and
only if $s$ is in $\Xi$. This measure is called the Riesz measure
with parameter $s$. \vskip0.1cm\n When $ s=(s_1,s_2,\ldots,s_r)\in
\Xi\setminus\displaystyle\prod_{i=1}^r\left]\displaystyle\frac{i-1}{2},+\infty\right[,$
the measure $R_s$ is concentrated on the boundary $\partial\Omega$
of $\Omega$ and when $s=(s_1,s_2,\ldots,s_r)$ is such that for all
$i$, $s_i>\displaystyle\frac{i-1}{2},$ the measure $R_s$ is
absolutely continuous with respect to the Lebesgue measure and is
given by
$$R_s(dx)=\frac{\Delta _{s-\frac{n}{r}}(x)}{\Gamma _{\Omega }(s)}
{\mathbf{1}}_{_{\Omega}}(x)(dx)$$ with
$n=\displaystyle\frac{r(r+1)}{2}$ the dimension of $E$
and
\begin{equation}\label{BAS1}
\Gamma_{\Omega }(s)=(2\pi
)^{\frac{r(r-1)}{4}}\displaystyle\prod_{i=1}^r\Gamma
\left(s_{i}-\frac{i-1}{2}\right).
\end{equation} It is shown
in Hassairi and Lajmi [6] that the measure $R_s$ generates a natural
exponential family if and only if $s_1\neq 0$. In this case
$$F=F(R_s)=\left\{R(s,\sigma )=\displaystyle\frac{e^{-<\sigma ,x>}}{\Delta _{s}(\sigma
^{-1})}R_s, \ \sigma\in\Omega, \ s\in\Xi, \ s_1\neq 0 \ \right\}.$$
The distribution $R(s,\sigma )$ is called the Riesz distribution
with parameters $s$ and $\sigma.$ It is shown in Hassairi and Lajmi
[6] that the Laplace transform of $R(s,\sigma)$ is defined for
$\theta$ in $\sigma-\Omega$, by
\begin{equation}\label{M9}
L_{R(s,\sigma)}(\theta)=
\displaystyle\frac{\Delta_s((\sigma-\theta)^{-1})}{\Delta_s(\sigma^{-1})}.
\end{equation}This implies that if $\sigma$ is an element of $\Omega$ and if $s$ and $s'$ are in $\Xi$,
then we have
$$R(s,\sigma)\ast R(s',\sigma)=R(s+s',\sigma).$$Let $\sigma$ be an element of $\Omega$. If $s$
satisfies the conditions $s_i>\displaystyle\frac{i-1}{2}$, for
$1\leq i\leq r,$ then the Riesz distribution is given by
$$
R(s,\sigma )(dx)=\displaystyle\frac{e^{-<\sigma ,x>}\Delta
_{s-\frac{n}{r}}(x)}{\Gamma _{\Omega }(s)\Delta _{s}(\sigma ^{-1})}
{\mathbf{1}}_{_{\Omega }}(x)(dx).$$ When $s_1=s_2=\ldots=s_r=p>0,$
$R(s,\sigma )$ reduces to the Wishart distribution with parameters
$p>\displaystyle\frac{r-1}{2}$ and $\sigma\in \Omega$,
$$W(p,\sigma)(dx)=\frac{1}{\Gamma _{\Omega}(p)\det (\sigma
^{-p})}e^{-<\sigma ,x>}\det(x)^{p-\frac{n}{r}}
{\mathbf{1}}_{_{\Omega }}(x)(dx),$$ with Laplace transform equal for
all $\theta\in \sigma-\Omega$, to
$$L_{W(p,\sigma)}(\theta)=\det\left(I_r-\sigma^{-1}\theta\right)^{-p}.$$
 \section{The mixture of
the Riesz distribution with respect to the multivariate Poisson}
Consider the Poisson distribution on $\natu^r$ with parameter
$\lambda=(\lambda_1,\lambda_2,\ldots,\lambda_r)\in (\reel^{+})^{r}$
\begin{equation}\label{FAT13}
\nu(dx)=e^{-\sum_{i=1}^r\lambda_i}\displaystyle\sum_{q\in\natu^r}\displaystyle\frac{\lambda^q
}{q!} \ \delta_q(dx),
\end{equation} where $q!=q_1! \ q_2!\ldots q_r!$ and $\lambda^q=\lambda_1^{^{q_1}}\lambda_2^{^{q_2}}
\ldots\lambda_r^{^{q_r}}.$ Then for all $\theta\in\reel^r,$
\begin{equation}\label{FATM5}
L_\nu(\theta)=\displaystyle\prod_{i=1}^r e^{\lambda_i
(e^{\theta_i}-1)}.
\end{equation} The variance function of the natural exponential
family generated by $\nu$ evaluated for $m\in ]0,+\infty[^{^{r}}$ is
equal to
$$V(m)=diag(m).$$ Note that this variance function is homogeneous of order
$1$ and consequently the natural exponential family $F(\nu)$ belongs
to the Tweedie scale of $\reel^r$ (see Hassairi and Louati
[7]).\vskip0.2cm\noindent
 Next, we will show that the mixture of the distribution $R(s,I_r)$ with
respect to $\nu$ is related to the modified Bessel function of the
first kind and of order $1$.\vskip0.1cm\n Consider the second-order
linear differential equation
\begin{equation}\label{NOUR5}
z^2y''(z)+zy'(z)+(z^2-\alpha^2)y(z)=0,
\end{equation} where $z$ is a complex variable
and $\alpha$ is a real parameter. One of the solution of
(\ref{NOUR5}) is the function $J(\alpha,z)$ known as the Bessel
function of the first kind of order $\alpha$, and defined for $z$ in
$\comp$ by the series
\begin{equation}\label{FATM4}
J(\alpha,z)=\displaystyle\sum_{k=0}^{+\infty}\displaystyle\frac{(-1)^k
\left(\displaystyle\frac{z}{2}\right)^{2k+\alpha}}{k!\Gamma(k+\alpha+1)}.
\end{equation}
In many applications, one frequently encounters two functions
$I(\alpha,z)$ and $K(\alpha,z)$ called respectively the modified
Bessel functions of the first and second kinds of order $\alpha$,
defined on the complex plane cut along the negative real axis by
\begin{equation}\label{FATM3}
I(\alpha,z)=i^{-\alpha}J(\alpha,iz),
\end{equation} and $$K(\alpha,z)=
\displaystyle\frac{\Pi}{2}\displaystyle\frac{I(-\alpha,x)-I(\alpha,x)}{\sin(
\alpha\Pi)}, \ \alpha\neq 0,\pm 1,\pm2,\ldots.$$ It is an immediate
fact that $I(\alpha,z)$ and $K(\alpha,z)$ are linearly independent
solutions of the differential equation:
\begin{equation}\label{FATM2}
z^2y''(z)+zy'(z)-(z^2+\alpha^2)y(z)=0.
\end{equation}
For more details about the Bessel functions, we refer to Lebedev
[10], page 108 and to Feller [4], Vol. II, page 58.\vskip0.1cm\n We
also need to establish the following result related to these
functions.
\begin{prop}\label{Proposition 1} Let $x>0$ and let $b>0,$ then
$$\displaystyle\sum_{k\in\natu}\displaystyle\frac{x^k}{k!
 \ \Gamma(k+b)}=x^{\frac{1-b}{2}} \ I(b-1,2\sqrt{x}).$$
\end{prop}

\begin{Pff} \
\vskip0.2cm\n Consider the function
\begin{equation}\label{NOUR2}
g_{_{b}}(x)=\displaystyle\sum_{k\in\natu}\displaystyle\frac{x^k}{k!
 \ \Gamma(k+b)}.
\end{equation}
Deriving $g_{_{b}}$ gives
$$g_{_{b}}'(x)=\displaystyle\sum_{k\in\natu}\displaystyle\frac{x^{k}}{k!
 \
 \Gamma(k+1+b)}=\displaystyle\sum_{k\in\natu}\displaystyle\frac{x^{k}}{k!(k+b)
 \ \Gamma(k+b)}.$$This implies that $$\left(x^bg_{_{b}}'(x)\right)'=\displaystyle\sum_{k\in\natu}\displaystyle\frac{x^{k+b-1}}{k!
 \ \Gamma(k+b)}=x^{b-1}g_{_{b}}(x).$$Consequently,
$g_{_{b}}(x)$ is a solution of the differential equation
$$x^b y''(x)+bx^{b-1}y'(x)-x^{b-1}y(x)=0.$$
Or equivalently, $g_{_{b}}(x)$ is solution of the equation
\begin{equation}\label{FAT24}
x y''(x)+by'(x)=y(x).
\end{equation}
Using now the fact that $I(\alpha,z)$ and $K(\alpha,z)$ are two
linearly independent solutions of (\ref{FATM2}), we deduce that
$x^{\frac{1-b}{2}} \ I(b-1,2\sqrt{x})$ and $x^{\frac{1-b}{2}} \
K(b-1,2\sqrt{x})$ are two linearly independent solutions of
(\ref{FAT24}). It follows that the general solution of (\ref{FAT24})
is of the form
$$c_1x^{\frac{1-b}{2}} \ I(b-1,2\sqrt{x})+c_2x^{\frac{1-b}{2}} \ K(b-1,2\sqrt{x}),$$
where $c_1$ and $c_2$ are two constants.\vskip0.1cm\n According to
(\ref{NOUR2}), (\ref{FATM4}) and (\ref{FATM3}), we deduce that
$$g_{_{b}}(1)=\displaystyle\sum_{k\in\natu}\displaystyle\frac{1}{k!
 \ \Gamma(k+b)}=I(b-1,2).$$Therefore $c_1=1$ and $c_2=0$.\vskip0.1cm\n Consequently, for
all $x>0,$
\begin{equation}\label{FAT30}
g_{_{b}}(x)=\displaystyle\sum_{k\in\natu}\displaystyle\frac{x^k}{k!
 \ \Gamma(k+b)!}=x^{\frac{1-b}{2}} \ I(b-1,2\sqrt{x}).
\end{equation}
\end{Pff}
\subsection{The probability density of the mixture of the Riesz distribution with respect to the
multivariate Poisson} Let
$\rho=\left(0,\frac{1}{2},\ldots,\frac{r-1}{2}\right)$,
$k=(k_1,k_2,\ldots,k_r)$ be in $\natu^r$ and let
$$\widetilde{R}_k=R(k+\rho,I_r).$$ Suppose that $k=(k_1,k_2,\ldots,k_r)$
has the multivariate Poisson distribution $\nu$ defined in
(\ref{FAT13}). For simplicity, in what follows, we will denote by
$\mu$ the mixture of $\widetilde{R}_k$ by $\nu$. The following
theorem gives the expression of $\mu$ in terms of the modified
Bessel function of the first kind.
\begin{theorem}\label{theorem0}
$$\mu(dx)=\displaystyle\frac{e^{-tr(x)}}{(2\pi)^{\frac{r(r-1)}{4}}\sqrt{\det(x)}}\displaystyle\prod_{i=1}^r
\displaystyle\frac{\sqrt{\lambda_i} \ e^{-\lambda_i}}{\sqrt
{\Delta_{i-1}(x)}} \
I\left(1,2\sqrt{\lambda_i\Delta_{e_i}(x)}\right){\mathbf{1}}_{_{\Omega}}(x)(dx),$$
where $\Delta_0(x)=1.$
\end{theorem}

\begin{Pff}
\vskip0.2cm\n Let $a=(a_1,a_2,\ldots,a_r)$ be an element of
$]0,+\infty[^{^{r}}$ and let
$$\widetilde{R}_{k,a}=R(k+\rho+a,I_r).$$ Denote
$\mu_{a}$ the mixture of $\widetilde{R}_{k,a}$ with respect to
$\nu$. Then
$$\mu_{a}(dx)=h_a(x){\mathbf{1}}_{_{\Omega}}(x)(dx),$$ where
\begin{equation}\label{NOUR3}
h_a(x)=e^{-\sum_{i=1}^r \lambda_i}\displaystyle\sum_{q\in \ \natu^{
r}}\displaystyle \frac{\lambda^q \
e^{-tr(x)}\Delta_{q+\rho+a-\frac{n}{r}}(x)}{q ! \
\Gamma_\Omega(q+\rho+a)},
\end{equation}
and
\begin{equation}\label{NOUR4}
\displaystyle\lim_{a\longrightarrow 0}
\mu_{a}(dx)=\mu(dx)=e^{-\sum_{i=1}^r
\lambda_i}\displaystyle\sum_{q\in \ \natu^{ r}}\displaystyle
\frac{\lambda^q \ e^{-tr(x)}\Delta_{q+\rho-\frac{n}{r}}(x)}{q ! \
\Gamma_\Omega(q+\rho)} \ {\mathbf{1}}_{_{\Omega}}(x)(dx).
\end{equation}
According to (\ref{FATM1}) and (\ref{BAS1}), (\ref{NOUR3}) may be
written as \vskip0.1cm\noindent
$h_a(x)=\displaystyle\frac{e^{-tr(x)}}{(2\Pi)^{\frac{r(r-1)}{4}}\sqrt{\det(x)}}\displaystyle\prod_{i=1}^r
\left(\Delta_i(x)^{-\frac{1}{2}}e^{-\lambda_i}\displaystyle\sum_{q_i\in\natu}\displaystyle\frac{\lambda_i^{q_i}}
{q_i! \ \Gamma(q_i+a_i)}\left(\displaystyle\frac{\Delta_i(x)}
{\Delta_{i-1}(x)}\right)^{q_i+a_i}\right).$\vskip0.1cm\hskip0.405cm
$=\displaystyle\frac{e^{-tr(x)}}{(2\Pi)^{\frac{r(r-1)}{4}}\sqrt{\det(x)}}\displaystyle\prod_{i=1}^r
\left(\Delta_i(x)^{-\frac{1}{2}}e^{-\lambda_i}\displaystyle\sum_{q_i\in\natu}\displaystyle\frac{\lambda_i^{q_i}}
{q_i! \ \Gamma(q_i+a_i)} \
\Delta_{e_i}(x)^{^{q_i+a_i}}\right).$\vskip0.1cm\hskip0.405cm
$=\displaystyle\frac{e^{-tr(x)}}{(2\Pi)^{\frac{r(r-1)}{4}}\sqrt{\det(x)}}\displaystyle\prod_{i=1}^r
\left(\Delta_i(x)^{-\frac{1}{2}}e^{-\lambda_i}
\Delta_{e_i}(x)^{^{a_i}}\displaystyle\sum_{q_i\in\natu}\displaystyle\frac{
\left(\lambda_i\Delta_{e_i}(x)\right)^{^{q_i}}} {q_i! \
\Gamma(q_i+a_i)}\right).$ \vskip0.1cm\noindent Therefore
\begin{equation}\label{FAT25}
h_a(x)=\displaystyle\frac{e^{-tr(x)}}{(2\Pi)^{\frac{r(r-1)}{4}}\sqrt{\det(x)}}\displaystyle\prod_{i=1}^r
\left(\displaystyle\frac{e^{-\lambda_i}} {\sqrt{\Delta_{i}(x)}}
 \ \Delta_{e_i}(x)^{^{a_i}}g_{a_i}\left(\lambda_i\Delta_{e_i}(x)\right)\right),
\end{equation} where $g_{a_i}$
is defined in (\ref{NOUR2}). Inserting now (\ref{FAT30}) in
(\ref{FAT25}), we obtain
$$h_a(x)=\displaystyle\frac{e^{-tr(x)}}{(2\Pi)^{\frac{r(r-1)}{4}}\sqrt{\det(x)}}\displaystyle\prod_{i=1}^r
\left(\displaystyle\frac{\lambda_i^{\frac{1-a_i}{2}}e^{-\lambda_i}}
{\sqrt{\Delta_{i}(x)}} \
\Delta_{e_i}(x)^{^{\frac{1+a_i}{2}}}I\left(a_i-1,2\sqrt{\lambda_i\Delta_{e_i}(x)
}\right)\right).$$ According to (\ref{NOUR4}), we deduce that
$$\mu(dx)=\displaystyle\frac{e^{-tr(x)}}{(2\Pi)^{\frac{r(r-1)}{4}}\sqrt{\det(x)}}\displaystyle\prod_{i=1}^r
\left(\displaystyle\frac{\sqrt{\lambda_i} \ e^{-\lambda_i}}
{\sqrt{\Delta_{i-1}(x)}} \
I\left(-1,2\sqrt{\lambda_i\Delta_{e_i}(x)}\right)\right){\mathbf{1}}_{_{\Omega}}(x)(dx).$$
Invoking the fact that for all $x>0,$ we have
$$I(1,x)=I(-1,x),$$(see Lebedev [10], page 110), the proof of the theorem is complete.
\end{Pff}
\subsection{The variance function of $F(\mu)$} In this subsection, we
study the natural exponential family $F$ generated by $\mu$. We
first give the Laplace transform of $\mu$, then we determine the
domain of the means and the variance function of the family
$F$.\vskip0.1cm\noindent In what follows, for $r\geq 1,$ we denote $
\kappa_{r}=\displaystyle\sum_{j=1}^r \displaystyle\frac{j}{2} \ e_j$
and we set $\kappa_{0}=0.$
\begin{theorem}\label{Theorem 1} For all
$\theta\in I_r-\Omega,$ we have
\begin{equation}\label{NOURM2}
L_{\mu}(\theta)=\Delta^\ast_{\kappa_{r-1}}(I_r-\theta)\exp
\left(\displaystyle\sum_{i=1}^r\lambda_i\left(
\Delta^\ast_{-e_{r-i+1}}(I_r-\theta)-1\right)\right).
\end{equation}
\end{theorem}

\begin{Pff} \
\vskip0.2cm\n Let $X_k$ be a random variable such that $X_k\sim
R(k+\rho,I_r).$ Then, according to (\ref{M9}), we have that for all
$\theta\in I_r-\Omega,$
$$L_{\mu}(\theta)=E\left(e^{\langle\theta,X_k\rangle}\right)=
E\left(E\left(e^{\langle\theta,X_k\rangle}\mid
k\right)\right)=E\left(\Delta_{k+\rho}((I_r-\theta)^{-1})\right).$$
Using (\ref{N1}) and (\ref{N0}), we can write
\vskip0.2cm\noindent$L_{\mu}(\theta)=E\left(
\Delta^\ast_{-(k+\rho)^\ast}(I_r-\theta)\right).$
\vskip0.2cm\hskip0.44cm$=E\left(\left(\Delta^{\ast
}_1(I_r-\theta)\right)^{k_{r-1}-k_r-\frac{1}{2}}\ldots\left(\Delta^{\ast
}_{r-1}(I_r-\theta)\right)^{k_1-k_2-\frac{1}{2}}\left(\Delta^{\ast
}_r(I_r-\theta)\right)^{-k_1}\right).$\vskip0.2cm\hskip0.44cm
$=\displaystyle\prod_{i=1}^{r-1}\left(\displaystyle\frac{1}
{\Delta^\ast_i(I_r-\theta)}\right)^{\frac{1}{2}}
E\left(\displaystyle\prod_{i=1}^{r}\left(\displaystyle\frac{\Delta^\ast
_{i-1}(I_r-\theta)}{\Delta^\ast
_{i}(I_r-\theta)}\right)^{k_{r-i+1}}\right).$\vskip0.2cm\hskip0.44cm
$=\displaystyle\prod_{i=1}^{r-1}\left(\displaystyle\frac{1}
{\Delta^\ast_i(I_r-\theta)}\right)^{\frac{1}{2}}
E\left(\displaystyle\prod_{i=1}^{r}\left(\Delta^\ast
_{-e_i}(I_r-\theta)\right)^{k_{r-i+1}}\right).$\vskip0.3cm\n It
follows that for all $\theta\in I_r-\Omega,$
\begin{equation}\label{FATM7}
L_{\mu}(\theta)=\displaystyle\prod_{i=1}^{r-1}\left(\displaystyle\frac{1}
{\Delta^\ast_i(I_r-\theta)}\right)^{\frac{1}{2}}
E\left(\displaystyle\prod_{i=1}^{r}e^{k_{r-i+1}\log(\Delta^\ast
_{-e_i}(I_r-\theta))}\right).
\end{equation}
Setting $\alpha(\theta)=\left(\log(\Delta^\ast
_{-e_r}(I_r-\theta)),\log(\Delta^\ast
_{-e_{r-1}}(I_r-\theta)),\ldots,\log(\Delta^\ast
_{-e_1}(I_r-\theta))\right)$, then (\ref{FATM7}) becomes
$$L_{\mu}(\theta)
=\displaystyle\prod_{i=1}^{r-1}\left(\displaystyle\frac{1}
{\Delta^\ast_i(I_r-\theta)}\right)^{\frac{1}{2}} E\left(e^{\langle
\alpha(\theta),k\rangle}\right).$$As $k$ has the multivariate
Poisson distribution $\nu,$ then
$$L_{\mu}(\theta)=\displaystyle\prod_{i=1}^{r-1}\left(\displaystyle\frac{1}
{\Delta^\ast_i(I_r-\theta)}\right)^{\frac{1}{2}} L_{\nu}
\left(\alpha(\theta)\right).$$ According to (\ref{FATM5}), we can
write for all $\theta\in I_r-\Omega,$
$$L_{\mu}(\theta)
=\displaystyle\prod_{i=1}^{r-1}\left(\displaystyle\frac{1}
{\Delta^\ast_i(I_r-\theta)}\right)^{\frac{1}{2}}
\displaystyle\prod_{i=1}^{r}e^{\lambda_i\left(e^{\log\left(\Delta^\ast
_{-e_{r-i+1}}(I_r-\theta)\right)}-1\right)}.$$Therefore
\begin{equation}\label{LM2}
L_{\mu}(\theta)
=\displaystyle\prod_{i=1}^{r-1}\left(\displaystyle\frac{1}
{\Delta^\ast_i(I_r-\theta)}\right)^{\frac{1}{2}}
\displaystyle\prod_{i=1}^{r}e^{\lambda_i\left(\left(\Delta^\ast
_{-e_{r-i+1}}(I_r-\theta)\right)-1\right)}.
\end{equation}On the other hand, using (\ref{N1}), we have $$\Delta^\ast_{\kappa_{r-1}}(I_r-\theta)=
\displaystyle\prod_{i=1}^{r-1}\left(\displaystyle\frac{1}
{\Delta^\ast_i(I_r-\theta)}\right)^{\frac{1}{2}}.$$ Inserting this
in (\ref{LM2}), we get (\ref{NOURM2}).
\end{Pff}

\begin{theorem}\label{theorem1} The domain of the means of the natural exponential family
$F=F(\mu)$ generated by the mixture $\mu$ is $\Omega.$
\end{theorem}

\begin{Pff} \
\vskip0.2cm\n From (\ref{NOURM2}), we deduce that for all $\theta\in
\Theta(\mu)=I_r-\Omega,$
\begin{equation}\label{NOUR8}
k_{\mu}(\theta)=\displaystyle\sum_{i=1}^r
\lambda_i\left(\left(\Delta^\ast
_{-e_{r-i+1}}(I_r-\theta)\right)-1\right)+\log(\Delta^\ast_{\kappa_{r-1}}(I_r-\theta)).
\end{equation}As for
all $i\in\{1,\ldots,r\}$, the map
$$\varphi_i: \ x\longmapsto
\log\Delta_i^\ast(x)$$ is differentiable on $\Omega$ and
\begin{equation}\label{MIN6}
\varphi'_i(x)=(P_i^\ast(x))^{-1},
\end{equation}then, for all $i\in\{1,\ldots,r\},$ we have
\begin{equation}\label{NOUR7}
\left(\Delta^\ast
_{-e_i}(x)\right)'=\left(\displaystyle\frac{\Delta^\ast
_{i-1}(x)}{\Delta^\ast _{i}(x)}\right)'=\Delta^\ast
_{-e_i}(x)\left((P_{i-1}^\ast(x))^{-1}-(P_i^\ast(x))^{-1}\right),
\end{equation}and for $r\geq 2,$ we have
\begin{equation}\label{NOURM1}
\left(\log(\Delta^\ast_{\kappa_{r-1}}(x))\right)'=
-\displaystyle\frac{1}{2}\displaystyle\sum_{i=1}^{r-1}(P_i^\ast(x))^{-1}
\end{equation}
Differentiating (\ref{NOUR8}) and taking into account (\ref{NOUR7})
and (\ref{NOURM1}), we get
\begin{equation}\label{M1}
k'_{\mu}(\theta)=
\displaystyle\sum_{i=1}^{r}\left(\lambda_{r-i+1}\Delta^\ast
_{-e_i}(I_r-\theta)-\lambda_{r-i}\Delta^\ast
_{-e_{i+1}}(I_r-\theta)+\displaystyle\frac{1}{2}\right)\left(P^\ast
_i(I_r-\theta)\right)^{-1} -\displaystyle\frac{1}{2}\left(P^\ast
_r(I_r-\theta)\right)^{-1},
\end{equation}where
$\lambda_0=0.$\vskip0.1cm\noindent Let $\theta\in I_r-\Omega,$ and
let $u$ be the unique element of ${\mathcal{T}}_{l}^+$ such that
$I_r-\theta=u^{\ast -1}(I_r).$ Then, for all $i\in\{1,\ldots,r\},$
we have
$$
\left(P^\ast_i(I_r-\theta)\right)^{-1}=\left(P^\ast_i(u^{\ast
-1}(I_r))\right)^{-1}=
\left(P^\ast_i\left((u(I_r))^{-1}\right)\right)^{-1}.$$
\vskip0.2cm\noindent According to (\ref{N15}), this  implies that
for all $i\in\{1,\ldots,r\},$
\begin{equation}\label{F2}
\left(P^\ast_i(I_r-\theta)\right)^{-1}=u\left(\displaystyle\sum_{j=r-i+1}^r
c_j\right).
\end{equation}
On the other hand, using (\ref{R1}), we can write for all
$i\in\{1,\ldots,r\},$
$$\Delta^\ast_{-e_i}(I_r-\theta)=\Delta^\ast_{-e_i}\left(u^{\ast
-1}(I_r)\right)=\Delta_{e_i^\ast}(u(I_r))=\Delta_{e_{r-i+1}}(u(I_r))
=\displaystyle\frac{\Delta_{r-i+1}(u(I_r))}{\Delta_{r-i}(u(I_r))}.$$
This with (\ref{Nour1}) imply that for all $i\in\{1,\ldots,r\},$
\begin{equation}\label{N2}
\Delta^\ast_{-e_i}=u_{r-i+1}^2,
\end{equation} where for all $i\in\{1,\ldots,r\}$, $u_i$ are defined in (\ref{III1}).\vskip0.1cm\noindent Using (\ref{F2}) and
(\ref{N2}), we deduce from (\ref{M1}) that
$$k'_{\mu}(\theta)=\displaystyle
\sum_{i=1}^{r}\left(\lambda_{r-i+1}u_{r-i+1}^2-\lambda_{r-i}u_{r-i}^2+
\displaystyle\frac{1}{2}\right) u\left(\displaystyle\sum_{j=r-i+1}^r
c_j\right)-\displaystyle\frac{1}{2} \
u\left(\sum_{i=1}^{r}c_i\right).$$ This after a standard
calculation, gives
\begin{equation}\label{M3}
k'_{\mu}(\theta)=\displaystyle\sum_{i=1}^{r}\left( \lambda_i
u_i^2+\displaystyle\frac{i-1}{2}\right)u(c_i)=u\left(\displaystyle\sum_{i=1}^{r}a_i(\theta)c_i\right),
\end{equation}
where
\begin{equation}\label{BAS2}
a_i(\theta)=\displaystyle\frac{i-1}{2}+\lambda_i u_i^2.
\end{equation} As the $a_i$ are strictly positive, we deduce
that
\begin{equation}\label{MI11}
k'_{\mu}(\Theta(\mu))=k'_{\mu}(I_r-\Omega)\subseteq\Omega.
\end{equation}
Conversely, consider $y\in \Omega,$ then using the Cholesky
decomposition, there exists a unique $w\in {\mathcal{T}}_{l}^+$ such
that\vskip0.1cm\noindent
$$y=w(I_r)=w\left(\displaystyle\sum_{i=1}^r
c_i\right)=w\left(
P\left(\displaystyle\sum_{i=1}^r\displaystyle\frac{1}{\sqrt{a_i(\theta)}}
c_i\right)\left(\displaystyle\sum_{i=1}^r a_i(\theta)
c_i\right)\right).$$ Let $\theta$ such that $I_r-\theta=u^{\ast
-1}\left(I_r\right)\in\Omega,$ where $u=w
\displaystyle\sum_{i=1}^r\displaystyle\frac{1}{\sqrt{a_i(\theta)}} \
c_i$. Then
$$y=u\left(\displaystyle\sum_{i=1}^r a_i(\theta) c_i\right).$$
This, using (\ref{M3}), gives $$y=k'_{\mu}(\theta)\in
k'_{\mu}(I_r-\Omega).$$ Therefore
$$\Omega\subseteq k'_{\mu}(I_r-\Omega)=k'_{\mu}(\Theta(\mu)).$$
This with (\ref{MI11}) imply that the domain of the means of the NEF
$F=F(\mu)$ is \vskip0.2cm$\hfill
M_{F}=k'_{\mu}(\Theta(\mu))=\Omega.$
\end{Pff}

\n The following theorem gives the variance function of the natural
exponential family $F=F(\mu).$
\begin{theorem}\label{theorem20}
For all $m\in \Omega,$\vskip0.1cm\noindent
$V_{F}(m)=-\displaystyle\frac{1}{2}P
\left(\displaystyle\sum_{i=1}^r\displaystyle\frac{1}{b_i(m)}
\left(\left(P^\ast_{r-i+1}(m^{-1})\right)^{-1}
-\left(P^\ast_{r-i}(m^{-1})\right)^{-1}\right)\right)$
\vskip0.1cm\hskip1.12cm $+\displaystyle\sum_{i=1}^{r}
\left(\displaystyle\frac{\lambda_{r-i+1} \ \Delta_{e_{r-i+1}}(m)}
{b_{r-i+1}(m)}-\displaystyle\frac{\lambda_{r-i} \
\Delta_{e_{r-i}}(m)}
{b_{r-i}(m)}+\displaystyle\frac{1}{2}\right)$\vskip0.1cm\hskip1.09cm$
\times\left[P\left(
\displaystyle\sum_{j=r-i+1}^r\displaystyle\frac{1}{b_j(m)}
\left(\left(P^\ast_{r-j+1}(m^{-1})\right)^{-1}
-\left(P^\ast_{r-j}(m^{-1})\right)^{-1}\right)\right)\right]$
\vskip0.1cm\hskip1.12cm$+\displaystyle\sum_{i=1}^{r}
\displaystyle\frac{\lambda_{r-i+1} \ \Delta_{e_{r-i+1}}(m)}
{(b_{r-i+1}(m))^3}\left[\left(\left(P^\ast_{i}(m^{-1})\right)^{-1}-
\left(P^\ast_{i-1}(m^{-1})\right)^{-1}\right)\right.$
\begin{equation}\label{MIN12}
\hskip4.945cm\left.\otimes
\left(\left(P^\ast_{i}(m^{-1})\right)^{-1}-
\left(P^\ast_{i-1}(m^{-1})\right)^{-1}\right)\right],
\end{equation}where for all $i\in\{1,\ldots,r\},$ $b_i(m)=
\displaystyle\frac{i-1}{4}+\sqrt{\left(\displaystyle\frac{i-1}{4}\right)^2+
\lambda_i\Delta_{e_i(m)}}.$
\end{theorem}
\vskip0.1cm\noindent Usually, for the calculation of the variance
function, we set $m=k'_{\mu}(\theta)$ and we determine its
reciprocal $\theta=\psi_\mu(m)$. This is difficult to do in the
present situation, however, we are able to determine
$\left(P^\ast_{i}(I_r-\psi_\mu(m))\right)^{-1}$,
$\Delta^\ast_{-e_i}(I_r-\psi_\mu(m)),$ and $a_i(\psi_\mu(m))$ where
$a_i$ is defined in (\ref{BAS2}).
\begin{theorem}\label{Theorem}
For all $i\in\{1,\ldots,r\}$, \vskip0.2cm\noindent$i) \
a_i(\psi_\mu(m))=
\displaystyle\frac{i-1}{4}+\sqrt{\left(\displaystyle\frac{i-1}{4}\right)^2+
\lambda_i\Delta_{e_i}(m)}.$\vskip0.2cm\noindent $ii) \
\left(P^\ast_i(I_r-\psi_\mu(m))\right)^{-1}=
\displaystyle\sum_{j=r-i+1}^{r}\displaystyle\frac{1}
{a_j(m)}\left[\left(P^\ast_{r-j+1}(m^{-1})\right)^{-1}
-\left(P^\ast_{r-j}(m^{-1})\right)^{-1}\right].$\vskip0.2cm\noindent
\begin{equation}\label{NA10}
\hskip-8.72cm iii) \ \Delta^\ast_{-e_i}(I_r-\psi_\mu(m))=
\displaystyle\frac{\Delta_{e_{r-i+1}}(m)}{a_{r-i+1}(m)}.
\end{equation}
\end{theorem}

\begin{Pff} \
\vskip0.2cm\n $i)$ As from Theorem \ref{theorem1},
$m=k'_{\mu}(\theta)$ is in $\Omega$, there exists a unique
$v\in{\mathcal{T}}_{l}^+$ such that $m=v(I_r).$ According to
(\ref{M3}), we
have$$v(I_r)=m=u\left(\displaystyle\sum_{i=1}^{r}a_i(\psi_\mu(m))
c_i\right)= u\left(
P\left(\displaystyle\sum_{i=1}^{r}\sqrt{a_i(\psi_\mu(m))} \
c_i\right) \left(\displaystyle\sum_{i=1}^{r}c_i\right)\right).$$
Therefore
$$v(I_r)= u\left(
P\left(\displaystyle\sum_{i=1}^{r}\sqrt{a_i(\psi_\mu(m))} \
c_i\right) \left(I_r\right)\right).$$ It follows that
$$v=u \ \displaystyle\sum_{i=1}^{r}\sqrt{a_i(\psi_\mu(m))} \ c_i.$$ Or equivalently,
\begin{equation}\label{NA3}
u=v \ \displaystyle\sum_{i=1}^{r}\displaystyle\frac{1}
{\sqrt{a_i(\psi_\mu(m))}} \ c_i.
\end{equation}
On the other hand, using (\ref{Nour1}), (\ref{BAS2}) becomes
$$a_i(\psi_\mu(m))=\displaystyle\frac{i-1}{2}+\lambda_i\Delta_{e_i}(u(I_r)).$$
Then using (\ref{NA3}), we can write
$$a_i(\psi_\mu(m))=
\displaystyle\frac{i-1}{2}+\lambda_i\Delta_{e_i}\left(\left(v \
\displaystyle\sum_{j=1}^r\displaystyle\frac{1}
{\sqrt{a_j(\psi_\mu(m))}} \
c_j\right)\left(\displaystyle\sum_{j=1}^r\displaystyle\frac{1}
{\sqrt{a_j(\psi_\mu(m))}} \ c_j  \
v^\ast\right)\right).$$\hskip2.336cm$= \displaystyle\frac{i-1}{2}+
\lambda_i\Delta_{e_i}\left(v \ \displaystyle\sum_{j=1}^r
\displaystyle\frac{1}{a_j(\psi_\mu(m))} \ c_j \
v^\ast\right).$\vskip0.2cm\hskip1.719cm$=
\displaystyle\frac{i-1}{2}+
\lambda_i\Delta_{e_i}\left(v\left(\displaystyle\sum_{j=1}^r
\displaystyle\frac{1}{a_j(\psi_\mu(m))} \
c_j\right)\right).$\vskip0.3cm\n Therefore $a_i(\psi_\mu(m))$
satisfies the equation
\begin{equation}\label{M6}
a_i(\psi_\mu(m))=\displaystyle\frac{i-1}{2}+\displaystyle\frac{\lambda_i
v_i^2}{a_i(\psi_\mu(m))},
\end{equation}
where $v_i$ is defined in (\ref{III1}).\vskip0.2cm\noindent As
$a_i(\psi_\mu(m))>0$, we deduce that
$$a_i(\psi_\mu(m))=\displaystyle\frac{i-1}{4}+\sqrt{\left(\displaystyle\frac{i-1}{4}\right)^2
+\lambda_i v_i^2}.$$ On the other hand, since $m=v(I_r),$ then using
(\ref{Nour1}), we have that
$$v_i^2=\Delta_{e_i}(m).$$
\vskip0.2cm\noindent Consequently, for all $i\in\{1,\ldots,r\},$
$$a_i(\psi_\mu(m))=
\displaystyle\frac{i-1}{4}+\sqrt{\left(\displaystyle\frac{i-1}{4}\right)^2+
\lambda_i\Delta_{e_i}(m)}.$$ $ii)$ With the notations used above, we
can write for all $i\in\{1,\ldots,r\},$\vskip0.2cm\noindent
$\left(P^\ast_i(I_r-\psi_\mu(m))\right)^{-1}=\left(P^\ast_i(u^{\ast
-1}(I_r))\right)^{-1}.$\vskip0.2cm\hskip2.84cm
$=u\left(\displaystyle\sum_{j=r-i+1}^r
c_{j}\right).$\vskip0.2cm\hskip2.84cm
$=\left(v\displaystyle\sum_{i=1}^r\displaystyle\frac{1}
{\sqrt{a_i(\psi_\mu(m))}} \ c_i\right)
\left(\displaystyle\sum_{j=r-i+1}^r c_j\right)
\left(\displaystyle\sum_{i=1}^r\displaystyle\frac{1}
{\sqrt{a_i(\psi_\mu(m))}} \ c_i
 \ v^\ast\right).$\vskip0.2cm\hskip2.84cm
$=v\displaystyle\sum_{j=r-i+1}^r\displaystyle\frac{1}{a_j(\psi_\mu(m))}
 \ c_j \ v^\ast.$\vskip0.2cm\hskip2.84cm
$=v\left(\displaystyle\sum_{j=r-i+1}^r\displaystyle\frac{1}{a_j(\psi_\mu(m))}
 \ c_j\right)$\vskip0.2cm\noindent Thus
\begin{equation}\label{N13}
\left(P^\ast_i(I_r-\psi_\mu(m))\right)^{-1}=\displaystyle\sum_{j=r-i+1}^r
\displaystyle\frac{1}{a_j(\psi_\mu(m))} \ v(c_j).
\end{equation}
\vskip0.2cm\noindent As $m=v(I_r)$, then for all
$j\in\{1,\ldots,r\},$ we have
$$v(c_j)=v\left(\displaystyle\sum_{i=j}^r c_i
-\displaystyle\sum_{i=j+1}^r
c_i\right)=\left(P^\ast_{r-j+1}(m^{-1})\right)^{-1}
-\left(P^\ast_{r-j}(m^{-1})\right)^{-1}.$$ Inserting this in
(\ref{N13}), we deduce that
\begin{equation}\label{MI18}
\left(P^\ast_i(I_r-\psi_\mu(m))\right)^{-1}=
\displaystyle\sum_{j=r-i+1}^{r}\displaystyle\frac{1}
{a_j(\psi_\mu(m))}\left[\left(P^\ast_{r-j+1}(m^{-1})\right)^{-1}
-\left(P^\ast_{r-j}(m^{-1})\right)^{-1}\right].
\end{equation}
Consequently
\begin{equation}\label{MI19}
\left(P^\ast_{i}(I_r-\psi_\mu(m))\right)^{-1}-\left(P^\ast_{i-1}(I_r-\psi_\mu(m))\right)^{-1}=
\displaystyle\frac{1}{a_{r-i+1}(\psi_\mu(m))}
\left(\left(P^\ast_{i}(m^{-1})\right)^{-1}-
\left(P^\ast_{i-1}(m^{-1})\right)^{-1}\right).
\end{equation}
$iii)$ According to (\ref{N1}) we have
$$\Delta^\ast_i(I_r-\psi_\mu(m))=\Delta^\ast_i(u^{\ast -1}(I_r))=\Delta^\ast_{\vartheta_i}(u^{\ast
-1}(I_r)),$$where $\vartheta_i=\displaystyle\sum_{j=1}^i e_j.$ Using
(\ref{R1}) and (\ref{NA3}), we deduce that \vskip0.2cm\noindent
$\Delta^\ast_i(I_r-\psi_\mu(m))=\Delta_{-\vartheta_i^\ast}(u(I_r)).$\vskip0.1cm\hskip2.16cm
$=\Delta_{-\vartheta_i^\ast}(uu^\ast).$\vskip0.1cm\hskip2.16cm$=\Delta_{-\vartheta_i^\ast}\left(\left(v
\displaystyle\sum_{j=1}^r\displaystyle\frac{1}
{\sqrt{a_j(\psi_\mu(m))}} \ c_j\right)
\left(\displaystyle\sum_{j=1}^r\displaystyle\frac{1}
{\sqrt{a_j(\psi_\mu(m))}} \ c_j \
v^\ast\right)\right).$\vskip0.1cm\hskip2.165cm$=\Delta_{-\vartheta_i^\ast}\left(v
\displaystyle\sum_{j=1}^r\displaystyle\frac{1} {a_j(\psi_\mu(m))} \
c_j \ v^\ast\right).$\\ Therefore
$$\Delta^\ast_i(I_r-\psi_\mu(m))=\Delta_{-\vartheta_i^\ast} \left(v\left(\displaystyle\sum
_{j=1}^r\displaystyle\frac{1}{a_j(\psi_\mu(m))} \
c_j\right)\right).$$ Thus, using (\ref{FATM1}), we obtain
$$\Delta^\ast_i(I_r-\psi_\mu(m))= \displaystyle\frac{\Delta_{r-i}
\left(v\left(\displaystyle\sum_{j=1}^r\displaystyle\frac{
1}{a_j(\psi_\mu(m))} \ c_j\right)\right)}{\Delta_{r}
\left(v\left(\displaystyle\sum_{j=1}^r\displaystyle\frac{
1}{a_j(\psi_\mu(m))} \ c_j\right)\right)}.$$It follows that for all
$i\in\{1,\ldots,r\},$
$$\Delta^\ast_{-e_i}(I_r-\psi_\mu(m))=\displaystyle\frac{\Delta^\ast_{i-1}(I_r-\psi_\mu(m))}{\Delta^\ast_{i}(I_r-\psi_\mu(m))}
=\displaystyle\frac{\Delta_{r-i+1}\left(v\left(\displaystyle\sum_{j=1}^r\displaystyle\frac{
1}{a_j(\psi_\mu(m))} \
c_j\right)\right)}{\Delta_{r-i}\left(v\left(\displaystyle\sum_{j=1}^r\displaystyle\frac{
1}{a_j(\psi_\mu(m))} \ c_j\right)\right)}.$$Using (\ref{Nour1}), we
deduce that for all $i\in\{1,\ldots,r\},$
$$\Delta^\ast_{-e_i}(I_r-\psi_\mu(m))=\displaystyle\frac{v_{r-i+1}^2}{a_{r-i+1}(\psi_\mu(m))}.$$
Therefore\vskip0.1cm
$$\Delta^\ast_{-e_i}(I_r-\psi_\mu(m))=
\displaystyle\frac
{\Delta_{r-i+1}(m)}{a_{r-i+1}(\psi_\mu(m))\Delta_{r-i}(m)}.$$\hskip6.8cm$=\displaystyle\frac
{\Delta_{e_{r-i+1}(m)}}{a_{r-i+1}(\psi_\mu(m))}.$
\end{Pff}

\vskip0.2cm\n We are now in position to give the variance function
of the natural exponential family $F$ stated in Theorem
\ref{theorem20}.\vskip0.2cm

\begin{Pff} \textbf{of Theorem \ref{theorem20}} \
\vskip0.2cm\n We have that for all $m\in M_F=\Omega,$
$$V_F(m)=k''_{\mu}(\psi_{\mu}(m)).$$ Differentiating (\ref{M1}) and using (\ref{NOUR7}) and the fact that
for all $i\in\{1,\ldots,r\}$ and $x\in \Omega,$
$$\left(\left(P^\ast_i(x)\right)^{-1}\right)'=-P\left(\left(P^\ast_i(x)\right)^{-1}\right),$$
we get for all $\theta\in I_r-\Omega,$ \vskip0.6cm\noindent$
k''_{\mu}(\theta)=-\displaystyle\frac{1}{2}P\left(\left(P^\ast_{r}
(I_r-\theta)\right)^{-1}\right)$
\vskip0.1cm\hskip0.79cm$+\displaystyle\sum_{i=1}^{r}\left(\lambda_{r-i+1}\Delta^\ast
_{-e_i}(I_r-\theta)-\lambda_{r-i}\Delta^\ast
_{-e_{i+1}}(I_r-\theta)+\displaystyle\frac{1}{2}\right)\left(P\left(\left(P^\ast_{i}
(I_r-\theta)\right)^{-1}\right)\right)$\vskip0.2cm \hskip0.79cm$
+\displaystyle\sum_{i=1}^{r}\lambda_{r-i+1}\Delta^\ast
_{-e_i}(I_r-\theta)\left(\left(P^\ast_{i}(I_r-\theta)\right)^{-1}-\left(P^\ast_{i-1}(I_r-\theta)\right)^{-1}\right)\otimes
\left(P^\ast_{i}(I_r-\theta)\right)^{-1}$\vskip0.2cm
\hskip0.79cm$-\displaystyle\sum_{i=1}^{r} \lambda_{r-i}\Delta^\ast
_{-e_{i+1}}(I_r-\theta)\left(\left(P^\ast_{i+1}(I_r-\theta)\right)^{-1}-\left(P^\ast_{i}(I_r-\theta)\right)^{-1}\right)
\otimes \left(P^\ast_{i}(I_r-\theta)\right)^{-1}.$ \vskip0.2cm
\noindent It follows that for all $\theta\in I_r-\Omega,$
\vskip0.2cm\noindent$
k''_{\mu}(\theta)=-\displaystyle\frac{1}{2}P\left(\left(P^\ast_{r}
(I_r-\theta)\right)^{-1}\right)$
\vskip0.1cm\hskip0.79cm$+\displaystyle\sum_{i=1}^{r}\left(\lambda_{r-i+1}\Delta^\ast
_{-e_i}(I_r-\theta)-\lambda_{r-i}\Delta^\ast
_{-e_{i+1}}(I_r-\theta)+\displaystyle\frac{1}{2}\right)\left(P\left(\left(P^\ast_{i}
(I_r-\theta)\right)^{-1}\right)\right)$\vskip0.3cm \hskip0.79cm$
+\displaystyle\sum_{i=1}^{r}\lambda_{r-i+1}\Delta^\ast
_{-e_i}(I_r-\theta)$ \vskip0.1cm\hskip0.79cm $ \times\left[
\left(\left(P^\ast_{i}(I_r-\theta)\right)^{-1}
-\left(P^\ast_{i-1}(I_r-\theta)\right)^{-1}\right)\otimes
\left(\left(P^\ast_{i}(I_r-\theta)\right)^{-1}-
\left(P^\ast_{i-1}(I_r-\theta)\right)^{-1}\right)\right].
$\vskip0.2cm\n We need only to replace $\theta$ by $\psi_{\mu}(m)$,
then insert (\ref{NA10}), (\ref{MI18}) and (\ref{MI19}) to get the
expression of the variance function of $F=F(\mu)$ given in
(\ref{MIN12}).
\end{Pff}

\vskip0.2cm\centerline{\large{\bf References}}\bigskip

\noindent [1] Aalen, O.O. Modelling heterogeneity in survival
analysis by the compound Poisson

distribution. Ann. Appl. Probab. 2 (1992), 951-972.

\noindent [2] Bhattacharya, S.K., Holla, M.S. On a discrete
distribution with special reference to the

theory of accident proneness. J. Amer. Statist. Asso. 60 (1965),
1060-1066.

\noindent [3] Faraut, J., Kor\'{a}nyi, A. Analysis on Symmetric
Cones. Oxford University Press, Oxford.

(1994).

\noindent [4] Feller, W. An Introduction to Probability Theory and
its Applications. Vol. I and II,

Second Edition. New York: Wiley. (1971).

\noindent [5] Gindikin, S.G. Analysis on homogeneous domains.
Russian Math. Surveys. 29 (1964),

1-89.

\noindent [6] Hassairi, A., Lajmi, S. Riesz exponential families on
symmetric cones. J. Theoret.

Probab. 14 (4) (2001), 927-948.

\noindent [7] Hassairi, A., Louati, M. Multivariate stable
exponential families and Tweedie scale, J.

Statist. Plann. Inference. 139 (2009), 143-158 .

\noindent [8] Johnson, L.L., Kemp, A.W. and Kotz.S. Univariate
discrete distributions. Third Edi-

tion. Hoboken, N.J.: Wiley. (2005).

\noindent [9] Karlis, D. and Meligkotsidou, L. Finite mixtures of
multivariate Poisson distributions

with application. J. Statist. Plann. Inference. 137 (2007),
1942-1960.

\noindent [10] Lebedev, N.N. Special functions and their
applications. Dover Publications, Inc. New

York. Translated and edited by Richard Silverman. (1972).

\noindent [11] Letac, G. Lectures on natural exponential families
and their variance functions. Mono-

grafias de Matem\'{a}tica. 50, IMPA, Rio de Janeiro. (1992).

\noindent [12] Perline, R. Mixed Poisson distributions tail
equivalent to their mixing distributions.

Statistics  Probability Letters 38 (1988), 229-233.

\noindent [13] Seshadri, V. The Inverse Gaussian Distribution. A
case study in exponential families.

Oxford University Press, Oxford. (1994).

\noindent [14] Y{\o}rgensen, B. The Theory of Dispersion Models.
First Edition. Chapman and Hall.

(1997).
\end{document}